\documentclass[12pt]{article} 
\usepackage{amsfonts,latexsym,amssymb,hyperref}

\begin{document}

\newtheorem{definition}{\bf{Definition}}[section]
\newtheorem{lemma}{\bf{Lemma}}[section]
\newtheorem{proposition}{\bf{Proposition}}[section]
\newtheorem{theorem}{\bf{Theorem}}[section]
\newtheorem{remark}{\sc{Remark}}[section]

\title{Existence of proper weak solutions  to the Navier-Stokes-Fourier system}

\author{Luisa Consiglieri\footnote{Independent Researcher, Portugal.
\href{http://sites.google.com/site/luisaconsiglieri}{http://sites.google.com/site/luisaconsiglieri}}}
\maketitle

\begin{abstract}
The existence of proper weak solutions of the Dirichlet-Cauchy problem 
constituted by the Navier-Stokes-Fourier system
which characterizes the incompressible
homogeneous Newtonian fluids under thermal effects is studied.
We call  {\em proper} weak solutions such weak solutions that verify some
 local energy inequalities in analogy with
the suitable weak solutions for the Navier-Stokes equations.
Finally, we deal with some regularity for the temperature.
\end{abstract}



\section{Introduction}

We study  the existence and the regularity 
of weak solutions for the Dirichlet-Cauchy problem
constituted by the Navier-Stokes-Fourier system
which is  one of the prominent mathematical problems for full thermodynamical
systems describing flows of incompressible fluids
(see for example
\cite{co00,ampa,naumann} and the references therein).
The presence of the temperature dependent viscosity into the momentum equation
 and the Joule effect 
into the energy equation are the main contributions on the nonlinear behavior
of the coupled system of partial differential equations. 
Other  problem which
arises from fluid thermomechanics is the Boussinesq approximation describing the dynamics in the planet mantle
(see \cite{mnrr} and the references therein). 
 In fact, it is a simpler
 coupled system which does not include the Joule effect.
For instance, in \cite{dd,gg,hi} even the viscosity is considered constant.

Although recently the $C(]0,T[;C^{0,\alpha}(\Omega))$- regularity for any weak solution
is found to the non-stationary Stokes system in 3D \cite{nw},
the non-uniqueness and smoothness of weak solutions are related with the non-uniqueness of
Leray-Hopf's  solutions to the first initial boundary value problem for the three-dimensional
Navier-Stokes equations.
For any class of weak solutions of the non-stationary N-S equations
in three-dimensional spaces, it is known that a weak solution cannot be used as a test function in the
weak variational formulation.
Indeed, 
in the N-S-F system, the viscosity is a temperature dependent function
and the temperature is a solution of a parabolic equation with $L^1$-data due to the existence of the Joule effect.
  By these reasons, in the work \cite{bcm,roub}
(and some references therein) the N-S-F system  is constituted by momentum
and total energy  equations in order to the $L^1$ dissipative term in the internal energy equation is  obtained
by approximative methods using the weakly sequential lower semicontinuity of the norm.
On the other hand,
 in the N-S system the study of the class of suitable weak solutions is
being of interest since \cite{ckn,lin,she77,she80}, if
 homogeneous Dirichlet boundary conditions
are enforced (see for instance \cite{gue07,gkt2006}).

In this work, 
we prove the existence of weak solutions to the problem under study and that among the weak solutions
at least the existence of one proper weak solution is guaranteed, i.e.
such weak solutions that verify some local energy inequalities.
Moreover, we show a higher integrability
 of the gradient of the velocity, using the reverse
H\"older inequality, which
is used in the proof of existence of weak solutions in order to deal with the $L^1$- dissipative
term in the energy equation. 
Some regularity for the temperature appears as a direct
consequence.
We refer to \cite{LT} for the proof of 
the higher
summability of the gradient of the velocity in the stationary case.
The regularity of strong two-dimensional solutions
 was solved by the author in the paper
\cite{c09} if some smallness on the data is considered. 
The study of the partial regularity for the velocity of the fluid
will be an ongoing future work.

Let  $\Omega \subset \mathbb R^n$ be a bounded open domain sufficiently regular and $T>0$.
Let us consider the boundary-value problem of the N-S-F system:
\begin{eqnarray}
\begin{array}c
\partial_t {\bf u}- {\rm div }(\mu(\theta) D{\bf u})+({\bf u}\cdot\nabla ){\bf u} ={\bf  f} -
\nabla p\quad\mbox{ in }
Q_T:=\Omega\times]0,T[;\\
{\rm div }~ {\bf u} = 0\quad\mbox{ in }
Q_T;
\end{array}
\label{U1}
\\
\partial_t \theta- {\rm div }( k(\theta) \nabla \theta)+{\bf u}\cdot
\nabla\theta =\mu(\theta)  |D{\bf u}|^2\quad\mbox{ in }
Q_T;
\label{U2}
\\
{\bf u}\big\vert_{t=0}= {\bf u}_0,\qquad
\theta\big\vert_{t=0}= \theta_0\quad \mbox{ in }\Omega;\label{U3a}\\
{\bf u}={\bf 0},\qquad
\theta= 0\quad\mbox{ on  }\partial \Omega\times ]0,T[,
\label{U3}
\end{eqnarray}
where $p$ denotes the pressure,  $\mu$ the viscosity,
$\theta$ the temperature, $\bf u$  the
velocity of the fluid and
$D{\bf u}={1\over2}(\nabla {\bf u}
 +\nabla{\bf u}^T)$, 
$\bf f$ denotes the given external body forces, $k$ denotes
the conductivity and it is assumed constant, $k(\theta)\equiv k$.
Notice that the assumption of a constant conductivity is not necessary
for the proof of existence of weak solutions.
It is a sufficient condition for the proof of existence of the local energy inequalities
to the temperature solution.
For simplicity, the constant density is assumed equal to one,
and we do not consider the existence of the external heat
source, since the heating dissipative term is the main mathematical difficulty.
The product of two tensors is given by $D:\tau=D_{ij}\tau_{ij}$ and the norm by $|D|^2=D:D$.
 
The initial conditions are given in (\ref{U3a}).
For the sake of clarity we found convenient that 
the
boundary conditions which are given in (\ref{U3}) 
are assumed homogeneous Dirichlet conditions.

The outline of the paper is as follows: in next section we establish
the appropriate functional framework and we present the main results.
The Section \ref{s3} is devoted to 
the proof of the existence result (Theorem \ref{th1}).
In Section \ref{r3}, we prove some regularity result (Proposition \ref{p1}).

\section{Assumptions and main results}

Here we assume that $\Omega \subset
 \mathbb R^n$ is a bounded open set such that its boundary
$\partial\Omega\in C^{2}$.
In the framework of Lebesgue and Sobolev spaces, for $1\leq q\leq \infty$,
we introduce \cite{La} 
$${\bf J}^{1,q}_0(\Omega)=\{{\bf u}\in{\bf W}^{1,q}_0(\Omega):\ \nabla\cdot{\bf u}=0
\mbox{ in }\Omega\}$$
with norm $$\Vert\cdot\Vert_{1,q,\Omega}=\Vert\nabla\cdot\Vert_{q,\Omega},$$
and we set 
\[W^{2,1}_q(Q_T):=\{v\in L^q(0,T;W^{2,q}(\Omega)):\ \partial_t v \in L^q(Q_T)\}
.\]
For any set $A$, we  write $(u,v)_A:=\int_A uv$ whenever
$u\in L^q(A)$ and $v\in L^{q'}(A)$, where $q'=q/(q-1)$ is the conjugate exponent to $q$,
or simply $(\cdot,\cdot)$ whenever there exists no confusing  at all, and we use
 the symbol
$\langle\cdot,\cdot\rangle$ to denote
 a generic duality pairing,
not distinguished between scalar and vector fields.
We denote by bold the vector spaces of vector-valued or tensor-valued functions. 
 
The following assertions on data are assumed as well as the following
assumptions on the physical parameters appearing in the equations are
established.
\begin{description}
\item[(A1)] ${\bf f}:Q_T\rightarrow\mathbb R^n$ is given such that
${\bf f}\in {\bf L}^{2(1+\epsilon_0)}(Q_T)$ with $\epsilon_0>0$;
\item[(A2)] 
$\mu: {\mathbb R}\to {\mathbb R}$ is a continuous function such that
\begin{equation}
0<\mu_\# \le \mu(s)\le \mu^\# ,\qquad \forall s\in {\mathbb R};
\label{defmu}
\end{equation}

\item[(A3)]
$ 
{\bf u}_0\in {\bf L}^{2}(\Omega),$ $ \theta_0\in  L^{1}(\Omega)
$ such that
\begin{equation}\label{defeo}
{\rm div }~ {\bf u}_0 = 0\quad\mbox{ in }\Omega;\quad{\rm ess}\inf_{x\in \Omega}\theta_0(x)\geq 0.
 \end{equation}
\end{description}

\begin{definition}
We say that the triple $({\bf u},p,\theta)$ is a weak solution to the Navier-Stokes-Fourier (N-S-F) problem 
 (\ref{U1})-(\ref{U3}) in $Q_T$ if 
\begin{eqnarray*}
{\bf u} \in {\mathcal U}:=L^{\infty}(0,T;{\bf L}^2(\Omega))\cap L^2(0,T;
{\bf J}^{1,2}_0(\Omega)),
 \quad p\in L^{(n+2)/n}(Q_T), \\
\theta \in  {\mathcal E}:=L^{\infty}(0,T;L^1(\Omega))\cap L^q(0,T;W^{1,q}_0 (\Omega)),\quad q<{n+2\over n+1},\\
\partial_t{\bf u} \in {\mathcal X}:=L^{2}(0,T;{\bf W}^ {-1,2 }(\Omega))\cap L^ {(n+2)/n }(0,T;{\bf W}^{-1,(n+2)/n}
(\Omega)),\\
\partial_t\theta \in  L^1(0,T;W^{-1,\ell} (\Omega)),\quad {1\over \ell}={n\over q( n+1)}+{n\over 2(n+2)}
\end{eqnarray*}
and satisfies the variational formulation
\begin{eqnarray}
\label{wmotion}
\langle \partial_t {\bf u},{\bf v}\rangle+\int _{Q_T} \Big(\mu(\theta
)D{\bf u}:D{\bf v}+({\bf u}\cdot\nabla ) {\bf u}\cdot{\bf v} \Big) =
\int _{Q_T} ( {\bf f}\cdot {\bf v}+ p\ {\rm div~}{\bf v}) \nonumber\\
\forall {\bf v} \in L^\infty(0,T;{\bf W}^{1,\infty}_0(\Omega)),\qquad
{\bf u}\big\vert_{t=0}={\bf u}_0\mbox{ in }\Omega;\\
\langle \partial_t \theta,\phi\rangle+\int _{Q_T} \Big( k \nabla \theta-\theta
{\bf  u}\Big)\cdot \nabla \phi dx dt
=\int _{Q_T} \mu(\theta)| D{\bf u}|^2  \phi dx dt\nonumber\\
\forall \phi \in L^{\infty}(0,T; W^{1,\infty}_0(\Omega)),\qquad
\theta\big\vert_{t=0}= \theta_0 \mbox{ in }\Omega.
\label{wheat}
\end{eqnarray}
\end{definition}
The embedding ${\mathcal X}\hookrightarrow {L}^{(n+2)/n}(0,T;{\bf W}^{-1,(n+2)/n}(\Omega))$ occurs for $n=2,3$.
For every ${\bf u}\in{\mathcal U}\hookrightarrow {\bf L}^{2(n+2)/n}({Q_T})$, the convective term verifies
$ ({\bf u}\cdot\nabla){\bf u}\in{\bf L}^{(n+2)/(n+1)}({Q_T})$ and consequently
$({\bf u}\cdot\nabla){\bf u}\cdot{\bf u}\not\in L^1({Q_T})$. Moreover,
the advection term ${\bf u}\cdot\nabla\theta\not\in L^1({Q_T})$ if 
$\theta\in
{\mathcal E}$ and $q< (n+2)/(n+1)$.
Remark that $\theta\in
{\mathcal E}\hookrightarrow L^ {q(n+1)/n}({Q_T})$ and $ \theta{\bf u}\in{\bf L}^\ell({Q_T})$
for $\ell\geq 1$, i.e. $q>1$ if $n=2$ and $q\geq 15/14 $ if $n=3$
($q\geq 2n(n+2)/((n+4)(n+1))$).

\begin{definition}
We say  that a weak solution $({\bf u},p,\theta)$ of the N-S-F  problem is proper in the following sense.
The local energy inequality holds
\begin{eqnarray}
\int _\Omega {| {\bf u}(x,t)-{\bf a}|^2\over 2}\varphi^2(x,t)dx+
\int _{Q_t=\Omega\times ]0,t[} \mu(\theta
)|D{\bf u}|^2 \varphi^2 dx\, d\tau\leq\nonumber\\
\leq 2\int_{Q_t}\mu(\theta)\varphi D{\bf u}:(({\bf u}-{\bf a})\otimes\nabla \varphi) 
+\int_{Q_t}{|{\bf u}-{\bf a}|^2}
\varphi(\partial_t\varphi+{\bf u} \cdot
\nabla\varphi)+\nonumber\\
+2\int _{Q_t}  p\varphi{\bf u}\cdot\nabla\varphi dx d\tau+\int_{Q_t}{\bf f}\cdot( {\bf u}-{\bf a})
\varphi^2 dx d\tau,\qquad
\label{e1}
\end{eqnarray}
for all ${\varphi} \in {C}^{\infty}_0(Q_T)$, a.e. $t\in ]0,T[$ and for any ${\bf a}\in\mathbb R^n$,
 and two more local energy inequalities hold
\begin{eqnarray}
\int _\Omega \left(\sqrt{\zeta+\theta^{2}}\psi\right)(x,t)dx+\zeta
k\int _{Q_t}{|\nabla\theta|^2\over (\zeta+\theta^2)^{3/2}} \psi
dxd\tau\leq\nonumber\\
\leq\int _ {Q_t} \sqrt{\zeta+\theta^{2}}
\left( \partial_t\psi+k\Delta\psi+
{\bf u}\cdot \nabla\psi
\right)+\int _{Q_t} \mu(\theta)| D{\bf u}|^2  {\theta\psi\over(\zeta+\theta^2)^{1/2}},
\label{e2}\\
\xi
k\int _{Q_t}{|\nabla\theta|^2\over (1+\theta)^{\xi+1}} \psi dx\,d\tau
+\int _{Q_t} \mu(\theta)| D{\bf u}|^2  {\psi\over(1+\theta)^{\xi}}dx\,d\tau\leq\nonumber\\
\leq
{1\over 1-\xi}\int _\Omega\left( {(1+\theta)^{1-\xi}}\psi\right)(x,t)dx-
\nonumber\\
-{1\over 1-\xi}\int _ {Q_t} {(1+\theta)^{1-\xi}}
\left( \partial_t\psi+k\Delta\psi+
{\bf u}\cdot \nabla\psi
\right)dxd\tau,\label{e3}
\end{eqnarray}
for any $\zeta>0$ and $0<\xi<1$, and
for all $\psi \in C^{\infty}_0(Q_T)$ such that $\psi\geq 0$, a.e $t\in ]0,T[$. 
\end{definition}

\begin{remark}
Two local energy inequalities are stated for the temperature
because from (\ref{e2}) if we take $\zeta\rightarrow 0^+$ we obtain 
$\theta \in L^{\infty}(0,T;L^1_{\rm loc}(\Omega))$ and from (\ref{e3})
after some computations we get
$\nabla\theta\in{\bf L}^q_{\rm loc} (Q_T)$.
\end{remark}

\begin{theorem}[$n=2,3$]\label{th1}
Under the assumptions (A1)-(A3),
 the N-S-F problem defined by (\ref{wmotion})-(\ref{wheat})
 has proper weak solutions.
Moreover, $\nabla{\bf u}\in {\bf L}^{2(1+\epsilon)}_{\rm loc}(Q_T)$ for $0<\epsilon< \min\{(4-n)/(3n),
1/(n+2),\epsilon_0\}$ and $\theta\geq 0$ in $Q_T$.
\end{theorem}

Some interior regularity is proved.
\begin{proposition}[Interior regularity]\label{p1}
Let  $({\bf u},p,\theta)$ be a weak solution in accordance to Theorem \ref{th1}.
Then $\theta\in W^{2,1}_{1+\epsilon,{\rm loc}}(Q_T)$.
In particular, $\theta\in {L}^{(n+2)(1+\epsilon)/(n-2\epsilon)}_{\rm loc}(Q_T)$.
\end{proposition}

Henceforth $C$  will denote different positive constants depending on the data,
but not on the unknown functions $\bf u$, $p$ or $\theta$.

\section{Proof of Theorem \ref{th1}}
\label{s3}

The proof of Theorem \ref{th1} is split in  each subsection.
In sections \ref{ss1}, \ref{ss2} and \ref{ss3} for reader's convenience we delineate the main arguments
concerning the existence of approximate solutions (for details, see \cite{bcm}).
The sections \ref{ss4}, \ref{ss6} and \ref{ss5} are new and
they are the main contributions for the desired
existence result.

\subsection{The Faedo-Galerkin method}
\label{ss1}

For $\nu,\varepsilon>0$ fixed,
there exists $\{{\bf u}^{N,M},p^{N,M},\theta^{N,M}\}_{N,M\in\mathbb N}$ being of the form
\cite[Carath\'eodory Theorem]{w}
\begin{eqnarray*}
&&{\bf u}^{N,M}\in\langle{\bf w}^1,\cdots,{\bf
w}^N\rangle\Leftrightarrow {\bf u}^{N,M}
(x,t)=\sum_{j=1}^Nc_j^{N,M}(t){\bf w}^j(x),\\
&&p^{N,M}
(x,t)={\mathcal F}_\varepsilon({\bf u}^{N,M}(x,t)),\\
&&\theta^{N,M}\in\langle w^1,\cdots,w^M\rangle\Leftrightarrow \theta
^{N,M}
(x,t)=\sum_{j=1}^Md_j^{N,M}(t){w}^j(x)
\end{eqnarray*}
where
 $\{({\bf
w}^j,w^j)\}_{j\in \mathbb N}$ is a basis of ${\bf W}^{1,\beta}_0(\Omega)\times
W^{1,\beta}_0(\Omega)$ with $\beta>n$,
$\mathcal F_\varepsilon$ is the  continuous functional such that maps
${\bf u}\in {\bf  W}^{1,2}_0(\Omega)$ into $p\in W^{2,2}(\Omega)$ which
is the solution of the  homogeneous Neumann problem for the Laplace equation (see, for instance, \cite{galdi})
\begin{eqnarray*}
\varepsilon\Delta p(t)={\rm div}\,{ \bf u}(t) &&\textrm{ in }\Omega\\
\nabla p(t)\cdot {\bf n}=0 &&\textrm{ on }\partial\Omega\\
\int_\Omega p(t) dx=0,&&
\end{eqnarray*}
which satisfies,  a.e. $t\in]0,T[$,
\begin{eqnarray*}
&&\varepsilon\|p(t)\|_{2,2,\Omega}\leq C(\Omega)\|{\bf u}(t)\|_{1,2,\Omega},
\quad\forall {\bf u}(t)\in{\bf W}^{1,2}_0(\Omega)\\
&&\varepsilon\|p(t)\|_{1,r,\Omega}\leq
C(\Omega,r)\|{\bf u}(t)\|_{r,\Omega},\quad\forall {\bf u}(t)\in{\bf W}^{1,2}_0(\Omega)\cap
{\bf L}^r(\Omega),\qquad r>1.
\end{eqnarray*}
The functions ${\bf c}^{N,M}=(c^{N,M}_1,\cdots,c^{N,M}_N)$
and $d^{N,M}=(d^{N,M}_1,\cdots,d^{N,M}_M)$
solve the following system of ordinary differential equations,
for every $M,N\in {\mathbb N}$,
\begin{eqnarray}\label{g1}
{d\over dt}({\bf u}^{N,M},{\bf w}^j)-(\mathcal M_\nu({\bf u}^{N,M})
\otimes {\bf u}^{N,M},\nabla{\bf w}^j)+(\mu(\theta^{N,M})D{\bf u}^{N,M},D{\bf w}^j)-\nonumber\\
-(\mathcal F_\varepsilon({\bf u}^{N,M}),\nabla \cdot{\bf w}^j) =({\bf f},{\bf w}^j),\quad
j=1,\cdots,N;\qquad\\
{d\over dt}(\theta^{N,M},w^j)-(\theta^{N,M}\mathcal M_\nu({\bf u}^{N,M}),\nabla
w^j)
+k(\nabla \theta^{N,M},\nabla w^j)=\nonumber\\
=(\mu(\theta^{N,M})|D{\bf u}^{N,M}|^2,w^j)
,\quad
j=1,\cdots,M,\qquad\label{g2}
\end{eqnarray}
under the initial conditions ${\bf u}_0^N,\ \theta_0^{N,M}$ given by the
projections of ${\bf u}_0$ and the mollification $\theta^N_0$of $\theta_0$
(after extending $\theta_0$ by zero outside $\Omega$),
respectively, onto linear hulls of the base's vectors. 
Finally, 
$\mathcal M_\nu$ is the divergenceless part of the Helmholtz-mollification decomposition
\[\mathcal M_\nu ({\bf u}):=(\chi{\bf u})*\omega-\nabla h\]
with  $\omega$ denoting a mollifier with support in a ball of
radii $\nu$,
\[\chi(x)=\left\{\begin{array}{ll}
0&\mbox{ if dist}(x,\partial\Omega)\leq 2\nu\\
1&\mbox{elsewhere},
\end{array}\right.\]
and $h$ is the Helmholtz-mollification decomposition, that is,
\begin{eqnarray*}
&&\Delta h={\rm div}[(\chi{\bf u})*\omega]
 \mbox{ in }\Omega\\
&&\nabla h\cdot {\bf n}=0 \mbox{ on }\partial\Omega\\
&&\int_\Omega h dx=0.
\end{eqnarray*}

Moreover, the standard estimates hold, independently of $M$, \cite{LSU,Li}
\begin{eqnarray}
\sup_{t\in[0,T]} \|{\bf u}^{N,M}(t)\|^2_{2,\Omega}+\mu_\#\|D{\bf
  u}^{N,M}\|^2_{2,Q_T}+\nonumber\\
+\varepsilon\|\nabla p^{N,M}\|_{2,Q_T}^2
\leq \|{\bf  u}_0\|^2_{2,\Omega}+C\|{\bf f}\|^{2}_{2,Q_T
};\label{cotau}\\
\sup_{t\in[0,T]}
\|{\theta}^{N,M}(t)\|^2_{2,\Omega}+k\|\nabla{\theta}^{N,M}\|^2_{2,Q_T}
\leq \|\theta^N_0\|^2_{2,\Omega}+C(N);\label{cotaeq}\\
\left\|{d\over dt}c^{N,M}\right\|_{L^2(0,T)}\leq C(N);\\
\|\partial_t \theta^{N,M}\|_{2,W^{-1,2}(\Omega)}\leq C(N,\nu).\label{cotaet}
\end{eqnarray}
Hence, the initial value problem (\ref{g1})-(\ref{g2}) has a global-in-time solution 
and passes to the limit as $M$ tends to infinity ($N$ fixed) by standard compactness arguments, resulting
\begin{eqnarray}\label{g1n}
{d\over dt}({\bf u}^{N},{\bf w}^j)-(\mathcal M_\nu({\bf u}^{N})
\otimes {\bf u}^{N},\nabla{\bf w}^j)+(\mu(\theta^{N})D{\bf u}^{N},D{\bf w}^j)-\nonumber\\
-(\mathcal F_\varepsilon({\bf u}^{N}),\nabla\cdot {\bf w}^j) =({\bf f},{\bf w}^j),\quad
 j=1,\cdots,N,\quad\mbox{a.e. }t\in ]0,T[;\\
\langle \partial_t\theta^{N},\phi\rangle -(\theta^{N}\mathcal M_\nu({\bf u}^{N}),\nabla
\phi)
+k(\nabla \theta^{N},\nabla \phi)=\nonumber\\
=(\mu(\theta^{N})|D{\bf u}^{N}|^2,\phi)
,\quad\forall  {\phi}\in L^2(0,T;{H}^1_0(\Omega)).\qquad\label{g2n}
\end{eqnarray}

Moreover, the minimum principle holds, i.e. $\theta^N\geq 0 $ a.e. in $Q_T$.

In order to pass to the limit on $N$, when $N$ tends to infinity,
the estimates (\ref{cotaeq})-(\ref{cotaet}) are no more valid. Thus,
we recall the additional estimates (for details, see \cite{bcm})
\begin{eqnarray}\label{cotaei}
\|\theta^N\|_{\infty,L^1(\Omega)}&\leq&\mu^\#\|
|D{\bf u}^N|^2\|_{1,Q_T}+T\|\theta_0\|_{1,\Omega}+|Q_T|/2;\\\label{cotaer}
\|\nabla \theta^N\|^q_{q,Q_T}&\leq& C\left(\mu^\#\| |D{\bf u}^N|^2\|_{1,Q_T}
+\|\theta_0\|_{1,\Omega}\right)\times\|
\theta^N\|^{q(2-q)/(2n)}_{\infty,L^1(\Omega)}\\
\|\partial_t{\bf u}^N\|^{2}_{2, {\bf W}^{-1,2}(\Omega)} &\leq
& C\left(\nu+{1\over\varepsilon}+1\right)\|{\bf
    u}^N\|^2_{2, {\bf W}^{1,2}_0(\Omega)}+C
\|{\bf f}\|^{2}_{2,Q_T};\label{udt}\\
\|\partial_t \theta^{N}\|_{1,W^{-1,q}(\Omega)}&\leq&\|\nabla \theta^{N}
\|_{q,Q_T}+C(\nu)\|\theta^{N}\|_{q,Q_T}
+\mu^\#\| |D{\bf u}^N|^2\|_{1,Q_T},\quad\label{edt}
\end{eqnarray}
for every exponent $1<q<2-n/(n+1)$ (cf. \cite{bago,bg}).
Moreover we have
 the strong convergences 
\begin{eqnarray}
\label{strongpi}
\nabla p^N\rightarrow \nabla p_\varepsilon\mbox{ in }{\bf L}^2(Q_T);
\\
\label{convtau}
\mu(\theta^N)|D{\bf u}^N|^2\rightarrow\mu(\theta_\varepsilon)|D{\bf u}_\varepsilon|^2\mbox{ in }L^1(Q_T).
\end{eqnarray}

Then, the initial value problem (\ref{g1n})-(\ref{g2n})  passes to the limit as $N$ tends to infinity 
($\varepsilon$ fixed), concluding the below quasi-compressible approximative problem.

\subsection{The quasi-compressible approximative problem}
\label{ss2}

From Section \ref{ss1},
for each $ \varepsilon>0$,
 there exists $({\bf u}_\varepsilon,p_\varepsilon,\theta_\varepsilon)$ in
$(L^{\infty}(0,T;{\bf L}^2(\Omega))\cap L^2(0,T;{\bf W}^{1,2}_0(\Omega)))
\times
L^2(0,T;W^{2,2}(\Omega))\times \mathcal E$ such that
$\partial_t{\bf u}_\varepsilon\in L^{2}(0,T; {\bf W }^{-1,2}(\Omega))$ and 
$\partial_t{\theta}_\varepsilon\in L^{1}(0,T; {W }^{-1,q}(\Omega))$,  and
it satisfies
\begin{eqnarray}
\langle\partial_t{\bf u}_\varepsilon,{\bf v}\rangle+\int_\Omega
\nabla{\bf u}_\varepsilon:{\bf
v}\otimes  {\mathcal M}_\nu({\bf u}_\varepsilon) dx+
\int_\Omega\mu(\theta_\varepsilon)D{\bf u}_\varepsilon:D{\bf v}
dx=\nonumber\\
 = \int_\Omega({\bf
f}-\nabla p_\varepsilon)\cdot{\bf v}
dx,\quad \forall{\bf v}\in  {\bf W}^{1,2}_0(\Omega),
\quad\mbox{  a.e. }t\in]0,T[;\qquad\label{pbueps}\\
\langle\partial_t \theta_\varepsilon,\phi\rangle+\int_{Q_T} {\mathcal
M}_\nu({\bf u}_\varepsilon) \cdot\nabla \theta_\varepsilon\phi dz+k\int_{Q_T}
\nabla
\theta_\varepsilon\cdot\nabla\phi dz=\nonumber\\
=\int_{Q_T}\mu(\theta_\varepsilon)|D{\bf u}_\varepsilon |^2
\phi dz,\quad
\forall\phi\in L^\infty(0,T;W^{1,q'}_0(\Omega));\quad\label{eeps}\\
\varepsilon\int_\Omega
\nabla p_\varepsilon\cdot\nabla\phi dx+\int_\Omega{\rm div}\
{\bf u}_\varepsilon \phi dx=0,\ \forall\phi\in W^{1,2}(\Omega);\label{pieps}\\
{\bf u}_\varepsilon(\cdot,0)={\bf u}_0,\quad
\theta_\varepsilon(\cdot,0)=\theta_0.\qquad\nonumber
\end{eqnarray}

In order to pass to the limit, when $\varepsilon$ tends to zero  ($\nu$ fixed),
the estimate (\ref{udt}) is no more valid. To estimate
$p_\varepsilon$ independently of $\varepsilon$ we choose ${\bf
v}=\nabla\eta_\varepsilon$ as a test function in (\ref{pbueps}),
where $\eta_\varepsilon$ is the solution of the following
homogeneous Neumann problem for the Laplace equation (see, for
instance, \cite{galdi})
\begin{eqnarray*}
&&\Delta \eta_\varepsilon(t)=
p_\varepsilon(t)-{1\over
|\Omega|}\int_\Omega
p_\varepsilon(t)
dx \mbox{ in }\Omega\\
&&\nabla \eta_\varepsilon(t)\cdot {\bf n}=0 \mbox{ on }\partial\Omega\\
&&\int_\Omega \eta_\varepsilon (t)dx=0,
\end{eqnarray*}
which satisfies
\begin{eqnarray}\label{cotaeta}
\|\eta_\varepsilon(t)\|_{2,2,\Omega}^{2}\leq C(\Omega)\|{p_\varepsilon}(t)\|_{2,{\Omega}}^{2},\quad
\mbox{  a.e. }t\in]0,T[.
\end{eqnarray}

After some technical computations, it results  (for details, see \cite{bcm})
\[\|p_\varepsilon\|_{{2},Q_T}^{2}\leq C(\nu),\]
and consequently
\begin{equation}\label{udtt}
\|\partial_t{\bf u}_\varepsilon\|^{2}_{2, {\bf W}^{-1,2}(\Omega)}\leq
C(\nu).
\end{equation}
Note that 
$\mu(\theta_\varepsilon)|D{\bf u}_\varepsilon|^2\rightarrow\mu(\theta_\nu)|D{\bf u}_\nu|^2$ holds in $L^1(Q_T)$ 
(compare to (\ref{convtau})) and  it is not required
 the strong convergence of the pressure 
$
\nabla p_\varepsilon\rightarrow \nabla p_\nu$ in ${\bf L}^2(Q_T)$
(compare to (\ref{strongpi})) since (\ref{pieps})  holds for $\phi=p_\varepsilon$.
Then, the initial value problem (\ref{pbueps})-(\ref{eeps})  passes to the limit as $\varepsilon$ tends to zero 
($\nu$ fixed), concluding the below regularized problem.

\subsection{The regularized problem}
\label{ss3}

For each $\nu\in\mathbb N$,
 there exists $({\bf u}_\nu,p_\nu,\theta_\nu)$ in
${\mathcal U}
\times
L^{2}(Q_T)\times \mathcal E$ such that
$\partial_t{\bf u}_\nu\in
L^{2}(0,T;{\bf W }^{-1,2}(\Omega))$
and $\partial_t{\theta}_\nu\in L^{1}(0,T; {W }^{-1,q}(\Omega))$,  and
 it satisfies
\begin{eqnarray}
\langle\partial_t{\bf u}_\nu,{\bf v}\rangle+\int_\Omega \nabla{\bf  u}_\nu:
{\bf v}\otimes {\mathcal M}_\nu({\bf u}_\nu)dx+
\int_\Omega\mu(\theta_\nu)D{\bf  u}_\nu:D{\bf v}dx
=\nonumber\\
  =\int_\Omega{\bf f}\cdot{\bf v}dx+\int_\Omega p_\nu {\rm div}\,{\bf v}dx ,
  \quad\forall{\bf v}\in 
{\bf W }^{1,2}_0(\Omega),\quad
\mbox{  a.e. }t\in]0,T[;\qquad\label{pbum}\\
\langle\partial_t \theta_\nu,\phi\rangle+\int_{Q_T} {\mathcal
M}_\nu({\bf u}_\nu) \cdot\nabla \theta_\nu\phi dxdt+k\int_{Q_T}\nabla
  \theta_\nu\cdot\nabla\phi dxdt=\nonumber\\
=\int_{Q_T}\mu(\theta_\nu)|D{\bf u}_\nu|^2 \phi dxdt,\quad
\forall\phi\in L^\infty(0,T; W^{1,q'}_0(\Omega));\qquad\label{em}\\
{\bf u}_\nu(\cdot,0)={\bf u}_0,\quad \theta_\nu(\cdot,0)=\theta_0 .\qquad\nonumber
\end{eqnarray}

Now in order to pass to the limit in (\ref{pbum})-(\ref{em}) as $\nu$ tends to infinity,
neither (\ref{edt}) nor (\ref{udtt}) are   valid. Following the argument of \cite{bcm},
we decompose the pressure $p_\nu$ into $p_\nu:=p_{\nu,1}+p_{\nu,2}$
such that the two particular pressures, $p_{\nu,1}$ and $p_{\nu,2}$,
belong to bounded sets of $L^{(n+2)/n}(Q_T)$ and $L^{2}(Q_T)$, respectively,
independently of $\nu$.
For each $t\in ]0,T[$, $p_{\nu,1}$ is the unique solution to the problem
\begin{eqnarray*}
\int_\Omega p_{\nu,1}(t) dx=0,\qquad
-\langle p_{\nu,1}(t),\Delta\phi\rangle=\langle{\bf u}_\nu\otimes
\mathcal M_\nu({\bf u}_\nu)(t),D\nabla\phi\rangle,
\end{eqnarray*}
for all $\phi\in W^{2,2}(\Omega)$ such that $\nabla\phi\cdot {\bf n }=0$
on $\partial\Omega$,
and define $p_{\nu,2}:=p_\nu-p_{\nu,1}$. Then $p_{\nu,2}$ solves at each time level
\begin{eqnarray*}
\langle{p}_{\nu,2},\Delta \phi\rangle=
\int_\Omega\mu(\theta_\nu)D{\bf  u}_\nu:D{\nabla\phi}dx
-\int_\Omega{\bf f}\cdot\nabla\phi dx,
\end{eqnarray*}
for all $\phi\in W^{2,2}(\Omega)$ such that $\nabla\phi\in  {\bf W }^{1,2}_0(\Omega)$,
and the following estimate holds
\begin{eqnarray*}
\|p_{\nu,1}\|_{(n+2)/n,Q_T}\leq C;\qquad
\|p_{\nu,2}\|_{2,Q_T}\leq C.
\end{eqnarray*}
Thus we conclude
 the following uniform estimates
\begin{eqnarray}\label{ptu}
&&\|\partial_t{\bf u}_\nu\|_{\mathcal X}\leq
 C\left(1+\|{\bf
    u}_\nu\|^2_{2, {\bf W }^{1,2}_0(\Omega)}+
\|{\bf u}_\nu\|^{2}_{2(n+2)/n,Q_T}\right);\\
&&\|\theta_\nu{\bf u}_\nu\|_{1,{\bf W}
^{-1,\ell}(\Omega)}\leq C;\qquad \|\partial_t \theta_\nu\|_{1,W^{-1,\ell}
(\Omega)}\leq C,\label{pte}
\end{eqnarray}
for
$1/\ell=n/[q(n+1)]+n/[2(n+2)]<1$ since $q<(n+2)/(n+1)$ and $n<4$.
Then under compactness arguments \cite{si},
$({\bf u},p,\theta)$ satisfies the limit variational formulation
(\ref{wmotion}). However, with the above estimates we only obtain 
\begin{eqnarray*}
\langle\partial_t \theta ,\phi\rangle-\int_{Q_T} \theta{\bf u}\cdot\nabla
\phi dxdt+k\int_{Q_T}\nabla \theta\cdot\nabla\phi dxdt
\geq (\mu(\theta)|D{\bf u}|^2,\phi)
, \end{eqnarray*}
 for all $\phi\in
  C^1(\bar Q_T)$ such that $\phi\geq 0$ \cite{bcm}.
  So we will need to prove an additional estimate.
  
\subsection{Higher integrability of $\nabla{\bf u}_\nu$}
\label{ss4}

First let us remark that (\ref{pbum})-(\ref{em}) may be rewritten as
\begin{eqnarray}
\langle \partial_t {\bf u}_\nu,{\bf v}\rangle+\int _\Omega \Big(\mu(\theta_\nu
)D{\bf u}_\nu:D{\bf v}+
{\bf v}\otimes  \mathcal M_\nu({\bf u}_\nu) : \nabla{\bf u}_\nu \Big) dx =\nonumber\\=
\int _\Omega( {\bf f}\cdot {\bf v}+ p_\nu {\rm div~}{\bf v})dx
,\quad\mbox{a.e. }t\in [0,T]\quad
\forall {\bf v} \in      { \bf W   }^{1,2}_0(\Omega),\label{pbau}
\\
{\bf u}\big\vert_{t=0}= {\bf u}_0\mbox{ in }\Omega;\nonumber\\
\langle\partial_t \theta_\nu,\phi\rangle+\int _\Omega \Big( k \nabla \theta_\nu-
\theta_\nu \mathcal M_\nu({\bf  u}_\nu) \Big)\cdot \nabla \phi dx
=\int _\Omega \mu(\theta_\nu)| D{\bf u}_\nu|^2  \phi dx,
\label{pbae}\\\mbox{a.e. }t\in [0,T]\quad\forall \phi \in W^{1,\infty}_0(\Omega),\qquad
\theta\big\vert_{t=0}= \theta_0 \mbox{ in }\Omega.\nonumber
\end{eqnarray}

Let us prove the required estimate.
\begin{theorem}\label{apriori}
Let $({\bf u}_\nu,p_\nu,\theta_\nu) \in
\mathcal U\times
L^{(n+2)/n}(Q_T)\times \mathcal E$ verify the system (\ref{pbum})-(\ref{em}) 
then $\nabla{\bf u}_\nu\in {\bf L}^{2(1+\epsilon)}_{\rm loc}(Q_T)$ for $0<\epsilon< \min\{
(4-n)/(3n),1/(n+2),\epsilon_0\}$
and the following estimate holds:
\begin{eqnarray}
\| \nabla{\bf u}_\nu \|_{2(1+\epsilon),Q(z_0,R)} \leq C \left( 
\| \nabla{\bf u}_\nu \|_{2,Q_{}(z_0,2R)}+\| {\bf u}_\nu \|_{2(n+2)/n,Q_{}(z_0,2R)}
\right.\nonumber \\
\left.+\| {\bf f}\|_{2(1+\epsilon_0),Q_{}(z_0,2R)} + \|p_\nu
\|_{(n+2)/n,Q_{}(z_0,2R)}^{1/2}\right), \label{ca1}
\end{eqnarray}
for any cylinder $Q_{}(z_0,2R):=B(x_0,2R)\times]t_0-(2R)^2,t_0[\subset \subset Q_T$.
\end{theorem}

{\bf Proof.}
We denote the points of the space-time cylinder by $z = (x, t)$
 and employ a shorthand
notation $ dz = dx dt$. 
We  write $ -\hspace*{-0.35cm}\int_{Q(z,R)}v:={1\over 
\omega_n R^{n+2}}\int_{Q(z,R)} v$ whenever
$v\in L^1(Q_T)$, where $\omega_n$ denotes the measure of the $n$-dimensional unit ball in $\mathbb R^n$.
For every $z_0=(x_0,t_0)\in  Q_T$ and
$ R>0$ small enough such that $Q(z_0,2R)\subset\subset Q_T$,
in order to prove the higher integrability of $\nabla{\bf u}_\nu$ (cf.  \cite[Gehring Lemma]{al,ge,st}),
it is sufficient to show the following reverse estimate 
\begin{eqnarray}
-\hspace*{-0.45cm}\int_{Q(z_0,R)}|\nabla{\bf u}_\nu|^2dz
&\leq&\delta-\hspace*{-0.45cm}\int_{Q(z_0,2R)}|\nabla{\bf u}_\nu|^2dz\nonumber\\
&&+{B_1\over R^{n+1}}-\hspace*{-0.45cm}\int_{Q(z_0,2R)}|{\bf u}_\nu|^{2}dz+{B_2\over R} 
-\hspace*{-0.45cm}\int_{Q(z_0,2R)} |{\bf u}_\nu|^3 dz\nonumber\\
&&+{B_3}-\hspace*{-0.45cm}\int_{Q(z_0,2R)}|{\bf f}|^2dz
+R-\hspace*{-0.45cm}\int_{Q(z_0,2R)}|p_\nu|^{(n+2)/n}dz
\label{ge}\end{eqnarray}
for some $\delta\in[0,1[$ and  positive constants $B_1,B_2,B_3$, independent of ${\bf u}_\nu,p_\nu$ and $\theta_\nu$,
considering that, for $n<4$,
\begin{eqnarray*}
{\bf u}_\nu\in {\bf L}^{2(n+2)/n}(Q_T),\qquad {2(n+2)\over n}>2;\\
|{\bf u}_\nu|^{3/2}\in { L}^{4(n+2)/(3n)}(Q_T),\qquad {4(n+2)\over 3n}>2;\\
{\bf f}\in {\bf L}^{2(1+\epsilon_0)}(Q_T),\qquad \epsilon_0>0;\\
|{p}_\nu|^{(n+2)/n}\in {L}^{1}(Q_T).
\end{eqnarray*}
Thus, we take $0<\epsilon< \min\{2/n,
(4-n)/(3n),1/(n+2),\epsilon_0\}=\min\{
(4-n)/(3n),1/(n+2),\epsilon_0\}$.

Adapting  the argument used in \cite{ark},
let $\varphi\in C^\infty_0(Q(z_0,2R))$ be a cut-off function such that
$\varphi\equiv 1\mbox{ in }Q(z_0,R)$ and $|\nabla\varphi|\leq
C/R,$ $|\partial_t\varphi|\leq C/R^2$ in $Q(z_0,2R).$
Choose ${\bf v}=\varphi^2{ \bf u}_\nu$ as a test function in (\ref{pbau}),
then
\begin{eqnarray*}
\int_{Q_t}{1\over 2}{d\over
  dt}(\varphi^2|{\bf u}_\nu|^2)-\int_{Q_t}\varphi|{\bf u}_\nu|^2\partial_t
\varphi+\int_{Q_t}\mu(\theta_\nu)\varphi^2|D{\bf u}_\nu|^2+\\
+2\int_{Q_t}\mu(\theta_\nu)\varphi D{\bf u}_\nu:({\bf u}_\nu
\otimes\nabla\varphi)=I+\int_{Q_t}\varphi^2
{\bf  f }\cdot{\bf u}_\nu+2\int_{Q_t} p_\nu\varphi{\bf u}_\nu\cdot\nabla\varphi
\end{eqnarray*} 
where $I$ corresponds to the convective term
\begin{eqnarray*}
I&=&\int_{Q_t} \mathcal M_\nu({\bf u}_\nu)\otimes {\bf u}_\nu:(\varphi^2\nabla
{\bf u}_\nu+2\varphi\nabla\varphi\otimes {\bf u}_\nu)dz\\
&=&\int_{Q_t}\left( \mathcal M_\nu({\bf u}_\nu)\cdot\nabla\varphi\right)
|{\bf u}_\nu|^2 \varphi dz.
\end{eqnarray*}
Applying  the assumption (\ref{defmu}) and H\"older and Young inequalities, it follows
\begin{eqnarray*}
\int_{\Omega}{1\over 2}{
}(\varphi^2|{\bf u}_\nu|^2)(t)dx+\mu_\#\int_{Q_t}\varphi^2|D{\bf u}_\nu|^2dxd\tau
\leq \int_{Q_t}\varphi|{\bf u}_\nu|^2|\partial_t
\varphi|+\\
+\delta\int_{Q_t}\varphi|\nabla {\bf u}_\nu|^2+
C(\delta,\mu^\#)\int_{Q_t}\varphi| {\bf u}_\nu|^2
|\nabla\varphi|^2+|I|+\\+{1\over 2}\int_{Q_t}\varphi^2
|{\bf  f }|^2+{1\over 2}\int_{Q_t}\varphi^2|{\bf u}_\nu|^2
+ {nR\over n+2}\int_{Q_t} |p_\nu|^{(n+2)/n}\varphi^{(n+2)/(2n)}+\\+{2\over n+2}{1\over R^{n/2}}
\int_{Q_t}\varphi^{(n+2)/4}|{\bf u}_\nu|^{(n+2)/2 }|\nabla\varphi|^{(n+2)/2}.
\end{eqnarray*} 

By Korn inequality the following estimate holds
$$\int_{Q_t} \varphi^2|\nabla {\bf u}_\nu|^2\leq 2\int_{Q_t }
\varphi^2|D{\bf u}_\nu|^2+4\int_{Q_t}|\nabla\varphi|^2|{\bf u}_\nu|^2.$$

According to the properties of $\varphi$
it arises
\begin{eqnarray*}
\mu_\#\int_{Q(z_0,R)}|\nabla{\bf u}_\nu|^2\leq 
\delta\int_{Q(z_0,2R)}|\nabla{\bf u}_\nu|^2+{C(\delta,\mu^\#,n)\over
  R^{n+1}}\int_{Q(z_0,2R)}|{\bf u}_\nu|^2dz+\\
+|I|+C\int_{Q(z_0,2R)}|{\bf f}|^2
  dz
+{nR\over n+2}\int_{Q(z_0,2R)}|p_\nu|^{(n+2)/n} dz,
\end{eqnarray*}
 observing that
for $R<1$ we have $R^2>R^{n+1}=R^{n/2+(n+2)/2}$.
Since
\begin{eqnarray*}
|I|\leq {C\over R}\int_{Q(z_0,2R)}|\mathcal M_\nu({\bf u}_\nu)|| 
{\bf u}_{\nu}|^2dz\leq\\\leq
{C\over R}\|\mathcal M_\nu({\bf u}_\nu)\|_{3,Q(z_0,2R)}\|{\bf u}_\nu\|^2_{3,Q(z_0,2R)}
\leq{C\over R}\|{\bf u}_{\nu}\|^3_{3,Q(z_0,2R)}
\end{eqnarray*}
then we conclude (\ref{ge}).

\subsection{Existence of weak solutions}
\label{ss6}

From estimates (\ref{cotau}), (\ref{cotaei})-(\ref{cotaer}),
(\ref{ptu})-(\ref{pte}) and (\ref{ca1}), independent on $\nu$, 
 we can extract a subsequence,
still denoted by $({\bf u}_\nu,p_\nu,\theta_\nu)$, verifying
(\ref{pbum})-(\ref{em}) and
\begin{eqnarray}
{\bf u}_\nu\rightharpoonup {\bf u}\quad &\mbox{weakly* in }&L^\infty(0,T;
 {\bf L}^2(\Omega));\nonumber\\
\nabla {\bf u}_\nu\rightharpoonup\nabla {\bf u}\quad& \mbox{weakly in }&{\bf L}^2(Q_T);\label{convp}\\
\partial_t{\bf u}_\nu\rightharpoonup \partial_t{\bf u}\quad &\mbox{weakly in }&\mathcal X
;\nonumber\\
{\bf u}_\nu\rightarrow {\bf u}\quad& \mbox{strongly in }&{\bf L}^m(Q_T),\
 \mbox{ for }1\leq m<2(n+2)/n;\label{umu}\\
{\theta}_\nu\rightharpoonup {\theta}\quad& \mbox{weakly in }&L^q(0,T;W^{1,q}_0(\Omega)),\ 
\mbox{ for }1<q<2-n/(n+1);\nonumber\\
{\theta}_\nu\rightarrow \theta\quad& \mbox{strongly in }&L^m(Q_T),\ 
\mbox{ for }1\leq m<q(n+1)/n;\label{emm}\\
p_\nu\rightharpoonup {p}\quad &\mbox{weakly in }&L^{(n+2)/n}(Q_T);\nonumber\\
p_{\nu,1}\rightarrow {p}_1\quad &\mbox{strongly in }&L^{m}(Q_T),\ 
\mbox{for }1\leq m<(n+2)/n;\label{pi1}\\
p_{\nu,2}\rightharpoonup {p}_2\quad &\mbox{weakly in }&L^{2}(Q_T).\nonumber
\end{eqnarray}

Since div $ {\bf u}_\nu =0$ and ${\bf u }_\nu|_{\partial \Omega\times ]0,T[}=0$, it follows from the definition and
properties of Helmoltz decomposition that 
\begin{eqnarray*}
\mathcal M_{\nu}( {\bf u}_\nu) \in L^{\infty}(0,T; {\bf L}^m(\Omega)) \textrm{ for all }m\in [1,\infty),\\
\mathcal M_\nu( {\bf u}_\nu)\rightarrow {\bf u} \textrm{ strongly in } {\bf L}^{m}(Q_T),
 \mbox{ for }1\leq m<2(n+2)/n.
\end{eqnarray*}

Then, the initial value problem (\ref{pbum})-(\ref{em})  passes to the limit as $\nu$ tends to infinity,
concluding the required problem (\ref{wmotion})-(\ref{wheat}).

\subsection{Local energy inequalities}
\label{ss5}

\subsubsection{Proof of the local energy inequality (\ref{e1})}

Let  $\varphi\in C^\infty_0(Q_T)$
and choose ${\bf v}=\varphi^2({ \bf u}_\nu-{\bf a})$ as a test function in (\ref{pbau}) for an
arbitrary ${\bf a}\in \mathbb R^n$, arguing as in  Section \ref{ss4}
then
\begin{eqnarray*}
{1\over 2}
\|\varphi({\bf u}_\nu-{\bf a })\|^2_{2,\Omega}(t)
+\int_{Q_t}\mu(\theta_\nu)\varphi^2|D{\bf u}_\nu|^2 dxd\tau=\\
=\int_{Q_t}\varphi|{\bf u}_\nu-{\bf a}|^2\partial_t
\varphi dxd\tau
-\int_{Q_t}\mu(\theta_\nu)\varphi D{\bf u}_\nu:({\bf u}_\nu-{\bf a})
\otimes\nabla\varphi dxd\tau+\\+
\int_{Q_t} \left(\mathcal M_\nu({\bf u}_\nu)\cdot\nabla\varphi\right)
{\bf u}_\nu\cdot ({\bf u}_\nu-{\bf a}) \varphi dxd\tau
\\+\int_{Q_t}\varphi^2
{\bf  f }\cdot({\bf u}_\nu-{\bf a})+2\int_{Q_t} p_\nu\varphi({\bf u}_\nu-{\bf a})\cdot\nabla\varphi dxd\tau.
\end{eqnarray*} 

Next, arguing as in  Section \ref{ss6}
we can pass to the limit as $\nu$ tends to infinity
concluding (\ref{e1}).

\subsubsection{Proof of the local energy inequality (\ref{e2})}
\label{se2}

Going back  to the solution
 $({\bf u}^{N},p^{N},\theta^{N})\in
(L^{\infty}(0,T;{\bf L}^2(\Omega))\cap L^2(0,T;{\bf W}^{1,2}_0(\Omega)))
\times
L^2(0,T;W^{2,2}(\Omega))\times \mathcal E$
of  the initial value problem (\ref{g1n})-(\ref{g2n})
we can choose $\phi=\theta^N(\zeta+(\theta^N)^2)^{-1/2}\psi\in  L^2(0,T;{ W}^{1,2}_0(\Omega)))
$, for $\zeta>0$ and   a non-negative function
$\psi\in {C}^\infty_0(Q_T)$, as a test function in (\ref{g2n}).
First,  we calculate separately the following terms
\begin{eqnarray*}
\langle\partial_t\theta^{N},\phi\rangle
=\int_\Omega\left(\sqrt{\zeta+(\theta^N)^2}
\psi\right)(x,t)dx
-\int_{Q_t}
\sqrt{\zeta+(\theta^N)^2}
\partial_t\psi(x,\tau)dxd\tau\\
(\mathcal M_\nu({\bf u}^{N}),\nabla
\theta^{N}\phi)=\int_{Q_t}\mathcal M_\nu({\bf u}^{N})\cdot\nabla\left(
\sqrt{\zeta+(\theta^N)^2}
\right)\psi dxd\tau\\
=- \int_{Q_t}
\sqrt{\zeta+(\theta^N)^2}
\mathcal M_\nu({\bf u}^{N})\cdot \nabla\psi dxd\tau\\
(\nabla \theta^{N},\nabla \phi)=\zeta\int
_{Q_t}{|\nabla \theta^N|^2\over (\zeta+(\theta^N)^2)^{3/2}}\psi dxd\tau+
\int_{Q_t}\nabla\left(\sqrt{\zeta+(\theta^N)^2}
\right)\cdot\nabla\psi dxd\tau\\
=\zeta\int
_{Q_t}{|\nabla \theta^N|^2\over (\zeta+(\theta^N)^2)^{3/2}}\psi  dxd\tau- 
\int_{Q_t}  \sqrt{\zeta+(\theta^N)^2}
\Delta\psi dxd\tau.
\end{eqnarray*}
Thus, we conclude 
\begin{eqnarray}
\int_\Omega\left(\sqrt{\zeta+(\theta^N)^2}
\psi\right)(x,t)
dx+
\zeta k\int
_{Q_t}{|\nabla \theta^N|^2\over (\zeta+(\theta^N)^2)^{3/2}}\psi dxd\tau=\nonumber\\
=\int_{Q_t}\sqrt{\zeta+(\theta^N)^2}
\left( \partial_t\psi+k\Delta\psi+
\mathcal M_\nu({\bf u}^{N})\cdot \nabla\psi
\right)dxd\tau+\nonumber\\
+\int_{Q_t}\mu(\theta^N)|D{\bf u}^N|^2 {\theta^N\psi\over \sqrt{\zeta+(\theta^N)^2}
}dxd\tau.\label{ee}
\end{eqnarray}
Since $\psi\in C^\infty_0(Q_T)$ such that $\psi\geq 0$ we can pass successively
to the limit as $N$ tends to infinity and $\varepsilon$ and $\nu$ tend to zero,
obtaining (\ref{e2}).
Indeed, still denoting by
$({\bf u}^N,p^N,\theta^N)$ the subsequence extracted from the solutions $({\bf u}^N,p^N,\theta^N)$ of
(\ref{g1n})-(\ref{g2n}) it verifies
\begin{eqnarray*}
{\bf u}^N\rightharpoonup {\bf u}_\varepsilon\quad &\mbox{weakly* in }&L^\infty(0,T;{\bf L}^2
(\Omega));\\
\nabla{\bf u}^N\rightharpoonup \nabla{\bf u}_\varepsilon\quad& \mbox{weakly in }&{\bf L}^2(Q_T);\\
\partial_t {\bf u}^N\rightharpoonup \partial_t{\bf u}_\varepsilon\quad 
&\mbox{weakly in }&L^{2}(0,T; {\bf W }^{-1,2}(\Omega));\\
{\bf u}^N\rightarrow {\bf u}_\varepsilon\quad& \mbox{strongly in }&{\bf L}^m(Q_T),\ 
\mbox{ for }1\leq m<2(n+2)/n;\\
{\theta}^N\rightharpoonup {\theta}_\varepsilon\quad& \mbox{weakly in }&L^q(0,T;W^{1,q}_0(\Omega)),\
 \mbox{ for }1<q<2-n/(n+1);\\
{\theta}^N\rightarrow \theta_\varepsilon\quad& \mbox{strongly in }&L^m(Q_T),\ \mbox{ for }1\leq m<q(n+1)/n;\\
p^N\rightharpoonup {p_\varepsilon}\quad &\mbox{weakly in }&L^2(0,T;W^{2,2}(\Omega));\\
\end{eqnarray*}
where $({\bf u}_\varepsilon,p_\varepsilon,\theta_\varepsilon)$ is a solution of
(\ref{pbueps})-(\ref{eeps}). 
Moreover, considering that $\psi\geq 0$,
\[
{\nabla\theta^N\over(\zeta+(\theta^N) ^2)^{3/4}
}\rightharpoonup {\nabla\theta_\varepsilon\over(\zeta+(\theta_\varepsilon) ^2)^{3/4}
}\quad \mbox{weakly in }{\bf L }^2(Q_T),
\] the lower semicontinuity property  of the $L^2$-norm
and \[
\mu(\theta^N){\theta^N\over\sqrt{\zeta+(\theta^N) ^2}
}\rightharpoonup \mu(\theta_\varepsilon){\theta_\varepsilon\over\sqrt{\zeta+(\theta_\varepsilon) ^2}
}\quad \mbox{*-weakly in }{\bf L }^\infty(Q_T),
\] and (\ref{convtau}) hold,
then the equality (\ref{ee}) delivers to
the local energy inequality (\ref{e2}) satisfied by $({\bf u}_\varepsilon,p_\varepsilon,\theta_\varepsilon)$.

Next,  still denoting by
$({\bf u}_\varepsilon,p_\varepsilon,\theta_\varepsilon)$ the subsequence extracted from the solutions
 $({\bf u}_\varepsilon,p_\varepsilon,\theta_\varepsilon)$ of
(\ref{pbueps})-(\ref{eeps}) it
verifies
\begin{eqnarray*}
{\bf u}_\varepsilon\rightharpoonup {\bf u}_\nu\quad &\mbox{weakly* in }&L^\infty(0,T;{\bf L }
^2(\Omega));\\
{\bf u}_\varepsilon\rightharpoonup {\bf u}_\nu\quad& \mbox{weakly in }&L^2(0,T; {\bf W }
^{1,2}_0(\Omega));\\
\partial_t{\bf u}_\varepsilon\rightharpoonup \partial_t{\bf u}_\nu\quad
 &\mbox{weakly in }&L^{2}(0,T;{\bf W}^{-1,2}(\Omega));\\
{\bf u}_\varepsilon\rightarrow {\bf u}_\nu\quad& \mbox{strongly in }&{\bf L}^m(Q_T),\ \mbox{ for }1\leq m<2(n+2)/n;\\
{\theta}_\varepsilon\rightharpoonup {\theta}_\nu\quad& \mbox{weakly in }&L^q(0,T;W^{1,q}_0(\Omega)
),\ \mbox{ for }1<q<2-n/(n+1);\\
{\theta}_\varepsilon\rightarrow \theta_\nu\quad& \mbox{strongly in }&L^m(Q_T),\ 
\mbox{ for }1\leq m<q(n+1)/n;\\
p_\varepsilon\rightharpoonup {p_\nu}\quad &\mbox{weakly in }&L^{2}(Q_T).
\end{eqnarray*}
where $({\bf u}_\nu,p_\nu,\theta_\nu)$ is a solution of
(\ref{pbum})-(\ref{em}). 
Analogously to the above argument the local energy inequality (\ref{e2}) arises for 
 $({\bf u}_\nu,p_\nu,\theta_\nu)$.

In conclusion, the convergences in Section \ref{ss6} imply that
 the limit $\theta$ verifies $\theta\geq 0$ in $Q_T$
and the local energy inequality (\ref{e2}) holds for the weak solution
 $({\bf u},p,\theta)$ such that verifies (\ref{e1}).

\subsubsection{Proof of the local energy inequality (\ref{e3})}

Proceeding as in Section \ref{se2}, we go back  to the existence of the solution
 $({\bf u}^{N},p^{N},\theta^{N})\in
(L^{\infty}(0,T;{\bf L}^2(\Omega))\cap L^2(0,T;{\bf W}^{1,2}_0(\Omega)))
\times
L^2(0,T;W^{2,2}(\Omega))\times \mathcal E$
of  the initial value problem (\ref{g1n})-(\ref{g2n}). Next
we can choose $\phi=(1+\theta^N)^{-\xi}\psi\in  L^2(0,T;{ W}^{1,2}_0(\Omega)))
$, for $0<\xi<1$ and   a non-negative function
$\psi\in {C}^\infty_0(Q_T)$, as a test function in (\ref{g2n}).
Now, calculating separately the following terms
\begin{eqnarray*}
(\partial_t\theta^{N},\phi)
&=&{1\over 1-\xi }\int_\Omega(1+\theta^N)^{1-\xi}
\psi(x,t)dx-\\&&
-{1\over 1-\xi }\int_{Q_t}(1+\theta^N)^{1-\xi}
\partial_t\psi(x,\tau)dxd\tau\\
(\mathcal M_\nu({\bf u}^{N}),\nabla
\theta^{N}\phi)
&=&- {1\over 1-\xi}\int_{Q_t} (1+
\theta^{N})^{1-\xi}\mathcal M_\nu({\bf u}^{N})\cdot  {\nabla\psi}dxd\tau\\
(\nabla \theta^{N},\nabla \phi)&=&-\xi\int
_{Q_t}{|\nabla \theta^N|^2\over (1+\theta^N)^{\xi+1}}\psi
- {1\over 1-\xi }
\int_{Q_t}  { (1+\theta^N)^{1-\xi}}\Delta\psi.
\end{eqnarray*}
Thus, we conclude 
\begin{eqnarray*}
\xi k\int
_{Q_t}{|\nabla \theta^N|^2\over (1+\theta^N)^{\xi+1}}\psi dx d\tau
+\int
_{Q_t}\mu(\theta^N)|D{\bf u}^N|^2 {\psi\over (1+\theta^N)^{\xi}}dxd\tau=\\
={1\over 1-\xi }\int_\Omega(1+\theta^N)^{1-\xi}{\psi(x,t)}dx-\\
-{1\over 1-\xi}\int_{Q_t}
{ (1+\theta^N)^{1-\xi}}\left( \partial_t\psi+k\Delta\psi+
\mathcal M_\nu({\bf u}^{N})\cdot \nabla\psi
\right)dxd\tau.
\end{eqnarray*}
Since $\psi\in C^\infty_0(Q_T)$ such that $\psi\geq 0$ we can pass successively
to the limit as $N$ tends to infinity and $\varepsilon$ and $\nu$ tend to zero,
obtaining (\ref{e3}).

\section{Regularity of $\theta$ (Proposition \ref{p1})}
\label{r3}

Let  $({\bf u},p,\theta)$ be a weak solution in accordance to Theorem \ref{th1},
that is, it satisfies in the sense of distributions
\begin{equation}\label{he}
\partial_t \theta- k\Delta\theta=-{\bf u}\cdot
\nabla\theta +\mu(\theta)  |D{\bf u}|^2\quad\mbox{ in }
Q_T.
\end{equation}

Thanks to Theorem \ref{th1}, we have
\begin{eqnarray*}
&&{\bf u} \in L^{\infty}(0,T;{\bf L}^2(\Omega))\cap L^{2(1+\epsilon)}(0,T;{\bf W}^{1,2(1+\epsilon)}_{\rm loc}(\Omega))
 \hookrightarrow {\bf L}^{4(n+2)(1+\epsilon)/n}_{\rm loc}(Q_T),\\
&&|\nabla
{\bf u}|^2 \in  L_{\rm loc}^{1+\epsilon}(Q_T).
  \end{eqnarray*}
Thus we get
\[{\bf u}\cdot\nabla \theta\in  L_{\rm loc}^{\varkappa}(Q_T),\qquad 
 {1\over \varkappa}={n\over 4(n+2)( 1+\epsilon)}+{1\over q}.
\]
Since $q<(n+2)/(n+1)$ it follows
\[1<\varkappa<{4(n+2)(1+\epsilon)\over 5n+4+4\epsilon(n+1)},\qquad\mbox{for every }
\epsilon>0.
\]
In particular, $\varkappa>1+\epsilon$ if $\epsilon<(4-n)/(4n+4)$.

Let us split the proof into two cases.
\begin{description}
\item[n=2]
Since we have $
{\bf u}\cdot\nabla \theta\in  L_{\rm loc}^{1+\epsilon}(Q_T)$ supposing $\epsilon<1/6$,
the classical theory for the heat equation (\ref{he}) leads
$\theta\in W^{2,1}_{1+\epsilon,{\rm loc}}(Q_T) 
$.
  
\item[n=3] 
The classical theory for the heat equation (\ref{he}) leads
$\theta\in W^{2,1}_{\varkappa,{\rm loc}}(Q_T) \hookrightarrow
 L^ {5\varkappa/(5-\varkappa)  }(0,T;W^{1,5\varkappa/(5-\varkappa)}(\Omega))
$.
Thus we get $
{\bf u}\cdot\nabla \theta\in  L_{\rm loc}^{1+\epsilon}(Q_T)$ supposing $\epsilon<1/6$,
and the bootstrap argument guarantees the required result.
\end{description}

Finally, the last assertion is valid from the embedding 
\[W^{2,1}_{1+\epsilon}(Q_T)\hookrightarrow L^{\varsigma}(Q_T),\qquad {2\over n+2}={1\over 1+\epsilon}-{1\over\varsigma}.
\]


\begin{thebibliography}{99}
\bibitem{al}
A. Arkhipova and O. Ladyzhenskaya, 
On inhomogeneous incompressible fluids and reverse H\"older inequalities.
{\em Annali della Scuola Normale Superiore di Pisa - Classe di Scienze,
 Ser. 4}, {\bf 25}: 1-2 (1997),  51-67.

\bibitem{ark}
{ A. Arkhipova and O. Ladyzhenskaya},
On a modification of Gehring's Lemma.
{\em Zapiski Nauchnykh Seminarov POMI} {\bf 259} (1999), 7-18;
English transl.:
{\em J. Math. Sci. (N.Y.)} {\bf 109} :5 (2002), 1805-1813.

\bibitem{bago}
{L. Boccardo}, { A. Dall'aglio}, {T. Gallouet} and
{L. Orsina}, Nonlinear
parabolic equations with measure data. {\em Journal of
Functional Analysis} {\bf 147} (1997), 237-258.
 
\bibitem{bg}
{L. Boccardo} and {T. Gallouet}, Non-linear elliptic and parabolic equations
involving measure. {\em Journal of Functional Analysis} {\bf 87} (1989), 149-169.

\bibitem{gg}
G. Boling and Y. Guangwei,
On the suitable weak solutions for the Cauchy problem
of the Boussinesq equations.
{\em Nonlinear Analysis: Theory,  Methods \& Applications}  {\bf 26} :8 (1996), 1367-1385.

\bibitem{bcm}
 M. Bul\'\i cek,  L. Consiglieri and  J. M\'alek,
Slip boundary effects on unsteady flows of incompressible viscous heat conducting fluids
 with a non-linear internal energy-temperature relationship.
Preprint  CMAFUL, Pre-2007-014.

\bibitem{ckn}
L. Caffarelli, R. Kohn and L. Niremberg,
Partial regularity of suitable weak solutions of the Navier-Stokes equations.
{\em Comm. Pure  Appl. Math.}  {\bf 35} :6 (1982), 771-831.

\bibitem{co00}
{L. Consiglieri}, Weak solutions for a class of non-Newtonian
fluids with energy transfer.  {\it J. Math. Fluid Mechanics} 
{\bf 2} (2000), 267-293.

\bibitem{ampa}
{L. Consiglieri},
Friction boundary conditions on thermal incompressible viscous flows.
{\em Annali di Matematica Pura ed Applicata}  {\bf 187} :4 (2008), 647-665.

\bibitem{c09}
{L. Consiglieri},
{Regularity for the Navier-Stokes-Fourier system}.
{\em Differential Equations \& Applications}  {\bf 1} :4 (2009), 583--604.



\bibitem{LT}
{ L. Consiglieri} and {T. Shilkin},
  Regularity to stationary weak solutions for  generalized Newtonian
fluids with energy transfer. {\em Zapiski Nauchnyh  Seminarov POMI
} {\bf 271} (2000), 122-150;
English transl.:
{\em J. Math. Sci. (N.Y.)} {\bf 115} (2003), 2771-2788.

\bibitem{dd}
D. Donatelli,
On the artificial compressibility method for the Navier Stokes Fourier system.
 arXiv:0807.3842v1 [math.AP] 

\bibitem{galdi}
{G.P. Galdi}, {\em An introduction to the mathematical theory of the
Navier-Stokes equations. Linearized steady problems}. Springer Tracts in
Natural Philosophy {\bf 38}, New York 1994.

\bibitem{ge}
{F.W. Gehring},
The $L^p$-integrability of the partial derivatives
of a quasi conformal mappings. {\it Acta Math.} {\bf 130}
(1973), 265-277.

\bibitem{gue07}
J.-L. Guermond,
Faedo-Galerkin weak solutions of the Navier-Stokes equations
with Dirichlet boundary conditions are suitable.
{\em J.  Math. Pures Appl.}  {\bf 88}  (2007), 87-106.

\bibitem{gkt2006}
S. Gustafson, K. Kang and T.-P. Tsai,
Regularity criteria for suitable weak solutions of the
Navier-Stokes equations near the boundary.
{\em J. Differential Equations}  {\bf 226}  (2006), 594-618.

\bibitem{hi}
T. Hishida,
 Existence and regularizing properties of solutions for the nonstationary convection problem.
{\em Funkcialaj Ekvacioj} {\bf 34} (1991), 449-474.

\bibitem{La} {O.A. Ladyzhenskaya},
 Mathematical problems in the dynamics of a
viscous incompressible fluid. 2nd rev. aug. ed., "Nauka", Moscow,
1970: English transl. of 1st ed.,{\em The mathematical theory of
viscous incompressible flow.} Gordon and Breach, New York 1969.

\bibitem{LSU}
{O.A. Ladyzhenskaya, V.A. Solonnikov and N.N. Uraltseva}, Linear
and quasilinear equations of parabolic type.
 Translations of Mathematical Monographs, Vol. 23 American Mathematical Society, Providence, R.I. 1967.

\bibitem{lin}
F. Lin,
A new proof of the Caffarelli-Kohn-Niremberg theorem.
{\em Comm. Pure  Appl. Math.}  {\bf 51} :3 (1998), 241-257.

\bibitem{Li}
{ J.L. Lions}, {\em Quelques m\'ethodes de r\'esolution des probl\`emes aux
limites non lin\'eaires.} Dunod et Gauthier-Villars, Paris 1969.

\bibitem{mnrr}
{ J. M\'alek, J. Ne\v cas, M. Rokyta}
 and {M. Ru\v zi\v cka}, {\em Weak
and Measure-valued solutions to evolutionary PDEs.} Chapman and Hall,
London 1996.


\bibitem{naumann}
J. Naumann, On the existence of weak solutions to the equations 
of non-stationary motion of heat-conducting incompressible viscous fluids.
{\em Math. Meth. Appl. Sci.} {\bf 29} (2006), 1883-1906.

\bibitem{nw} 
J. Naumann and J. Wolf,
On the interior regularity of weak solutions
to the non-stationary Stokes system.
{\em J. Glob. Optim.} {\bf 40} (2008), 277-288.
 
 \bibitem{roub}
T. Roub\'\i \v cek, On non-Newtonian fluids with energy transfer.
 {\em J.~Math.~Fluid~Mech.}
{\bf 11} (2009), 110-125 (published on line in vol.9, no.1 (June 2007), 1-16.). 

\bibitem{she77}
V. Scheffer,
Hausdorff measure and the Navier-Stokes equations.
{\em Comm.  Math. Phys.}  {\bf 55} :2 (1977), 97-112.


\bibitem{she80}
V. Scheffer,
The Navier-Stokes equations on a bounded domain.
{\em Comm.  Math. Phys.}  {\bf 73}  (1980), 1-42.

\bibitem{st}
E.W. Stredulinsky,
Higher integrability from reverse H\"older inequalities.
 {\em Indiana Univ. Math. J.} {\bf 29} (1980), 408-417.

\bibitem{si}
{J. Simon},  Compact sets in the space. {\em Annali Mat. Pura Appl.IV} {\bf
146} (1987), 65-96.

\bibitem{w}
{W. Walter,} {\em Ordinary differential equations}. Springer-Verlag,
 New York, Inc. 1998.

\end{thebibliography}
\end{document}